\documentclass[10 pt,leqno]{amsart}

\usepackage {amsfonts}
\usepackage{amsthm}
\usepackage{amssymb}
\usepackage{latexsym}
\usepackage{amsmath}
\usepackage{mathrsfs}
\usepackage{pifont}

\theoremstyle{plain}
\newtheorem{tw}{Theorem}[section]

\newtheorem {lem} [tw]{Lemma}
\newtheorem {prop}[tw] {Proposition}

\newtheorem{cor}[tw]{Corollary}

\theoremstyle{definition}
\newtheorem {deft}[tw] {Definition}
\newtheorem {rem} [tw]{Remark}

\newcommand{\bn}{\Bbb N}

\newcommand{\bz}{\Bbb Z}
\newcommand{\ld}{\ldots}

\newcommand{\algen} { A}

\newcommand{\com} {\mathcal C}
\newcommand{\hyp} {\mathcal D}
\newcommand {\Tr} {{\textrm{Tr}}}
\newcommand {\hte} {{\textrm{ht}}}

\newcommand{\Zr}{\bz_+^r}

\newcommand{\jr}{\overline{j}}
\newcommand{\kr}{\overline{k}}

\newcommand{\maxshape}{n}

\newcommand{\la}{\langle}
\newcommand{\ra}{\rangle}

\newcommand{\ot}{\otimes}

\newcommand{\wt}{\widetilde}


\begin{document}
\author{Adam Skalski}
\author{Joachim Zacharias}
\footnote{\emph{Permanent address of the first named author:}
Department of Mathematics, University of \L\'{o}d\'{z}, ul.
Banacha 22, 90-238 \L\'{o}d\'{z}, Poland.

AS acknowledges the support of the Polish KBN Research Grant 2P03A 03024.}
\address{School of Mathematical Sciences,  University of Nottingham,
Nottingham, NG7 2RD}
\email{pmxags@maths.nottingham.ac.uk}
\email{joachim.zacharias@nottingham.ac.uk  }

\title{\bf Entropy of shifts on higher-rank graph $C^*$-algebras}
\maketitle
\begin{abstract}
\noindent
Let $O_{\Lambda}$ be a higher rank graph $C^*$-algebra. For every $p \in \Zr$ there is a canonical completely positive map
$\Phi^p$ on $O_{\Lambda}$ and a subshift $T^p$ on the path space $X=\Lambda^\infty$. We show that $ht(\Phi^p)=h(T^p)$,
where $ht$ is Voiculescu's approximation entropy and $h$ the classical topological entropy.
For a higher rank Cuntz-Krieger algebra $O_{M}$ we obtain $ht(\Phi^p)= \log \textup{r}(M_1^{p_1}M_2^{p_2} \ldots M_r^{p_r})$,
$\textup{r}$ being the spectral radius. This  generalizes Boca and Goldstein's result for Cuntz-Krieger algebras.
\end{abstract}

\vspace*{1cm}
Higher rank graph $C^*$-algebras were introduced in \cite{kupa} as a generalization of higher
rank Cuntz-Krieger algebras defined in \cite{Steg}. Given a row finite rank-$r$ graph  $\Lambda$
with no sources
Kumjian and Pask \cite{kupa} define an infinite path space $\Lambda^{\infty}$ together
with a semigroup of continuous shift maps $T^p : \Lambda^{\infty} \to \Lambda^{\infty}$,
where $p \in \bz_+^r$. If, as we will always assume, the set of objects (alphabet)
$B=\Lambda^0$ is finite then $\Lambda^{\infty}$ is a zero-dimensional compact space
(\cite{kupa}). When $r=1$ and $\Lambda$ is the graph associated to the 0-1-matrix $A$,
$\Lambda^{\infty}$ is a subshift of finite type $X_A$ and $T_A=T^1$ is the classical
shift map on $X_A$.
Generally $\Lambda^{\infty}$ is a higher rank
subshift of finite type. In the higher rank Cuntz-Krieger case $\Lambda$ is given by a tuple of
0-1-matrices $M=(M_1 , \ldots , M_r)$ satisfying the conditions (H0)-(H3) of \cite{Steg}.
Note that following the tradition, we use the name `shift' (or `subshift')
both for the transformation of a dynamical system and for a system itself.

Now let $O_{\Lambda}$ be the rank-$r$ graph $C^*$-algebra associated to $\Lambda$ (\cite{kupa}).
There is a semigroup of canonical unital completely positive maps
$\Phi^p : O_{\Lambda} \to O_{\Lambda}$ which on $C(\Lambda^{\infty}) \subseteq O_{\Lambda}$
is given by $f \mapsto f \circ T^p$ (and therefore is often called a \emph{noncommutative shift action}).
In the one-dimensional Cuntz-Krieger case $\Phi_A (x)= \sum_{i=1}^n s_i xs_i^*$. Goldstein and Boca proved
that $ht(\Phi_A) = h(T_A)$ (\cite{Boca}), where  $ht$ is Voiculescu's approximation entropy
and $h(T_A)$ the classical topological entropy known to be equal to $\log \textrm{r}(A)$, where
$\textup{r}$ denotes the spectral radius.

In this note we  extend this result to higher rank graph $C^*$-algebras.  The techniques used
are similar to those in \cite{Boca} and \cite{Kerr}. The main difference lies in the fact
that there are many commuting shifts in different directions and we can ask about
\begin{enumerate}
\item the entropy of a shift $\Phi^p$, where $p \in \Zr$;
\item the entropy of the action of the whole semigroup $\Phi$ by shifts.
\end{enumerate}
In both cases the entropy is equal to the entropy of the corresponding
classical subshift. For higher-rank Cuntz-Krieger algebras we obtain the value
$\hte (\Phi^p) = \log \textrm{r}(M_1^{p_1}M_2^{p_2} \ldots M_r^{p_r})=:\log \textrm{r}(M^{p})$
in case 1 whereas in case 2 the entropy is always 0.

The plan of the note is as follows: in the first part we fix notations
and compute classical topological entropies of the transformations $T^p$ on rank-$r$ subshifts
of finite type. The entropies may be known, but we are not aware
of any reference and therefore a short proof is provided.
Then we proceed to recall the definition of
higher rank graph $C^*$-algebras $O_{\Lambda}$ and canonical shift maps $\Phi^p$ acting on them.
In section 2 we show that the entropy of the canonical shift
coincides with the entropy of the classical subshift. The result is established for the more general
topological pressure. In the third section we define the
topological (approximation) entropy for $\Zr$ actions and show that $ht(\Phi)=0$
whenever $r\geq 2$.

\section{Notations and basic facts}

\subsection{Rank-$r$ subshifts}

We follow the notation in \cite{Steg} and \cite{kupa}.
Throughout this note we fix $r \in \bn$ and a finite alphabet $B=\Lambda^0$.
For any $m,n \in \Zr$, where $m \leq n$ in the standard
partial ordering of $\Zr$, let $[m,n] := \{l \in \Zr: m\leq l \leq n\}$.
$e_1,\ldots,e_r$ denotes the standard basis of $\Zr$ and $e=e_1+\ld+e_r$.
${\bf N}^r$ denotes the  category with one object and morphisms $\Zr$.
A rank-$r$ graph $\Lambda$ is a small category with set of objects $\Lambda^0$ and shape functor $\sigma : \Lambda \to {\bf N}^r$
satisfying the factorisation property. The morphisms in $\Lambda$ may be thought of as (multidimensional) words.
Given $m \in \bz_+^r$ the set
$$
\Lambda_m=\{w \in \Lambda \mid \sigma(w) = m \} \;\;\;\;\;\; \text{ with cardinality } w_m = | \Lambda_m |
$$
consists of the morphisms or words of shape $m$.
The factorisation property says
that every word $w \in \Lambda_{m+n}$ is a unique product $w=uv$ of
a word $u \in \Lambda_m$ and $v \in \Lambda_n$, where $t(u)=o(v)$. Here $o$ and $t$ are the origin
and terminal maps also called range and source.
$\Lambda^0$ can thus be identified with words of shape 0
and for any word $w \in \Lambda_l$ and $k\in \Zr, k\leq l$ there are well-defined restrictions
$w|_{k]}$  and  $w|_{[k}$  with shape $k$ and $l-k$ respectively, given by $w = w|_{k]}w|_{[k}$.
For any $m \in \Zr$ and $\lambda \in \Lambda_0$ we write
\[
\Lambda_m(\lambda)=\{w \in \Lambda_m : o(w) = \lambda\}. \]
Note that the factorisation property has strong implications concerning the number of
words of a given shape, which in turn facilitate the entropy computations in what follows (see for
example formula \eqref{zero0} in the proof of Proposition \ref{coment}).

We assume that the rank-$r$ graph $\Lambda$ is \emph{row finite} (that is for any $m \in \Zr$ and $\lambda \in \Lambda_0$
the set $\Lambda_m(\lambda)$ is finite) and has \emph{no sources} (that is for any $m \in \Zr$ and $\lambda \in \Lambda_0$
the set $\Lambda_m(\lambda)$ is nonempty).
The infinite path space is defined by
$$
\Lambda^{\infty}= \big\{ \big( w_{(m,n)} \big)_{m \leq n} :  w_{(m,n)} \in \Lambda_{n-m} , \; w_{(m,n)}=w_{(m,k)} w_{(k,n)}
\forall_{m \leq k \leq n} \big\},
$$
slightly more concrete than in \cite{kupa}.
The origin $o(w)$ is defined for every
$w \in \Lambda^{\infty}$ and we have a composition with any $v \in \Lambda$ provided $t(v)=o(w)$.
As shown in \cite{kupa} $X = \Lambda^{\infty}$ with the topology generated by
cylinder sets
\[
Z_u = \{x \in X :x_{(0,\sigma(u))} = u\}, \;\;\; u \in W
\]
is a compact space. We can define a metric on $X$ inducing this topology as follows. For $j \in \bz_+ $ let $\jr= (j,\ldots,j) \in \Zr$ and
define $d(x,x)= 0$ if $x \in X$ and for all $x,y \in X, x \neq y$
\[
d(x,y) = \frac{1}{k+1}, \;\; \textrm {where}\; k =\min \{j\in \bz_+:  x|_{\jr]}
\neq y|_{\jr]}\}.
\]
so that $X$ is a compact zero-dimensional metric space. For any $p \in \Zr$ define a (continuous) shift $T^p:X\to X$ by
\[
T^p (x)_{(m,n)} = x_{(m+p,n+p)}.
\]
The prime examples of rank-$r$ graphs are the  Robertson-Steger graphs defined as follows.
Let $M_1,\ldots,M_r:B \times B \to \{0,1\}$ be matrices satisfying conditions
(H0)-(H3) of \cite{Steg}. Let $W_0=B$,
\[
W_m = \left\{w:[0,m] \to B :
\forall_{l\in [0,m]} \forall_{j=1,\ldots,r} \, l +e_j \in [0,m] \Rightarrow
M_j(w(l), w(l+e_j)) =1 \right\},
\]
$W= \bigcup_{m \in \Zr}W_m$ and $\sigma(w) = m $ if $w \in W_m$, where
$o(w)=w(0)$, $t(w)=w(\sigma(w))$. Then $W$ is a rank $r$-graph.
The factorisation property follows from the crucial condition (H1). The infinite path space $W^{\infty}$ is given by
\[
X_M = \left\{ x:\Zr \to B : \forall_{m\in \Zr} \forall_{j=1,\ldots,r} \, M_j(x(m), x(m+e_j)) =1
\right\}.
\]
\begin{deft}
We refer to the compact space $(X,d)$ equipped with the action of $\Zr$ by shifts $T$ as a rank-$r$ subshift and to
$(X_M,d)$ with corresponding $T$ as a rank-$r$ subshift of finite type
(given by $M=(M_1,\ldots,M_r)$).
\end{deft}
For subshifts of finite type we have the following explicit formula for the entropy of $T^p$.
\begin{prop} \label{coment}
The topological entropy of $T^p$ for a rank-$r$ subshift of finite type given by
 $M=(M_1, \ldots ,M_r)$ is equal to $ \log \textrm{r}(M^p) := \log \textrm{r}(M_1^{p_1}\cdot
 \cdots \cdot M_r^{p_r})$ and the topological entropy
of the action of $\Zr$ on any rank-$r$ subshift $X$ by $T$ is equal to $0$ (when $r \geq 2$).
\end{prop}

\begin{proof}
 We will use the definition of entropy given by R.\,E.\,Bowen, based on the notion of separating
subsets (\cite{Walt}, Section 7.2). If $n\in \bn$ define for $x,y \in X$
\[
d_n(x,y) = \max \{d(T^{lp}(x), T^{lp}(y)): 0\leq l\leq n\}.
\]
Put (for each $k\in \bn$) $\varepsilon_k = \frac{1}{k+2}$. Then
\[
d_n(x,y) > \varepsilon_k \Leftrightarrow \exists_{l\in\{0,\ldots,n\}} T^{lp} (x)|_{\kr]} \neq
T^{lp} (y)|_{\kr]},
\]
and in case $\bar{k} \geq p$ this is equivalent to
\[
x|_{\bar{k}+np]} \neq y|_{\bar{k}+np]}
\]
(this follows from the factorisation property). It is therefore easy to see that maximal $(n, \varepsilon_k)$-separating subsets of $X_M$ are
exactly sets of words whose restrictions to $[0,\bar{k}+np]$ are mutually different and exhaust
all possibilities. Therefore the cardinality of a maximal $(n, \varepsilon_k)$-separating
subset is exactly equal to $w_{\bar{k}+np}=|W_{\bar{k}+np}|$.
Observe that in general, for $l \in \Zr$,
\begin{eqnarray*}
w_l &=&
\textrm{card} \Big\{ w:[0,l]\to B: M_1 (w(0), w(e_1))=1, M_1 (w(e_1), w(2e_1))=1, \ldots \\
& &
 \text{\hspace{1.5cm}} \ldots ,M_2 (w(l_1 e_1), w(l_1e_1+e_2))=1,
\ldots , M_r( w(l-e_r), w(l))=1 \Big\} \\
&=&
\sum_{j_0 \in B}\sum_{j_1^1,\ldots,j_{l_1}^1\in B} \ld  \sum_{j_1^r,\ldots,j_{l_r}^r\in B}
M_1(j_0,j_1^1)M_1(j_1^1,j_2^1) \cdots M_1(j_{l_1-1}^1,j_{l_1}^1) \\
& & \text{\hspace{3cm}} M_2(j_{l_1}^1,j_1^2) \cdots M_2(j_{l_2-1}^2,j_{l_2}^2) \cdots M_r(j_{l_r-1}^r,j_{l_r}^r)\\
&& \\
&=&
\sum_{j_0,j_{l_r}^r\in B} M_1^{l_1} \cdots  M_r^{l_r} (j_0,j_{l_r}^r) =
\la e, M_1^{l_1} \cdots  M_r^{l_r} e\ra= \la e, M^{l} e\ra,
\end{eqnarray*}
where in the first equality property (H1) was used.
Further note that for all $l,m \in \Zr$
\[ \|M^l\| \leq \la e, M^l e \ra = w_l \leq w_{l+m} = \la e, M^{l+m} e \ra \leq \|M^{l+m}\| \cdot \|e\|_2^2
 \leq \|M^{l}\| \cdot \|M^m \| \cdot |\Lambda_0| \]
(the second inequality follows from the factorisation property and the assumption that $\Lambda$ has no sources).
As a consequence,
\begin{eqnarray}
\limsup_{n\in \bn} \frac{1}{n}\log s(n, \varepsilon_k)
&=&
\limsup_{n\in \bn} \frac{1}{n}\log w_{\bar{k}+np} =
\limsup_{n\in \bn} \frac{1}{n}\log \la e, M^{\bar{k}+np} e\ra \\
&=& \limsup_{n\in \bn} \frac{1}{n}\log \|M^{np}\| = \log \textrm{r}(M^p), \nonumber \label{logent}
\end{eqnarray}
for $k$ sufficiently large.
Thus
\[
h_{\textrm{top}} (T) = \sup_{\varepsilon >0} \left(\limsup_n \Big( \frac{1}{n}\log s(n,\varepsilon) \Big) \right) =\log \textrm{r}(M^p).
\]
For the second part observe that the maximal separating subsets of $X$ described above
are also maximal $(n,\varepsilon_k)$ separating sets for the action $T$ of $\Zr$ by shifts
(for defining entropy of $\Zr$-actions we use standard F\o lner sets $C_n=\{a\in \Zr: \max\{a_1,\ld,a_r\} \leq n\}$, cf.\:section 3).
This immediately yields (for $r\geq 2$)
\begin{equation} \label{zero0}
h_{\textrm{top}} (\alpha, \varepsilon_k) = \limsup_{n\in \bn} \frac{1}{n^r}\log w_{\overline{k+n}} \leq \limsup_{n\in \bn} \frac{1}{n^r}
\log |B|^{2(n+k)r}= 0,
\end{equation}
and therefore the entropy in question vanishes.
\end{proof}

The description of maximal $(n,\varepsilon_k)$-separating subsets is also valid for a general subshift
and we obtain the following formula for the topological pressure of $T$.
\begin{rem} \label{compre}
If $T$ is the canonical shift on $X=\Lambda^{\infty}$, $f\in C(X)_+$ and $k,n \in \bn$ then
\[
p_{\textrm{top}} (f,T^p,n,\varepsilon_k) = \sum_{u\in \Lambda_{\bar{k}+np}} \exp(\max\{f(x)+f(T^p(x))+\ld +
f(T^{np}(x)): x \in Z_u\}),
\]
whenever $p \in \Zr$ is such that $\bar{k} \geq p$.
\end{rem}
For various equivalent definitions of topological entropy and pressure of a continuous transformation of a compact space (or
for actions of amenable semigroups on a compact space) we refer to \cite{Walt} (or \cite{Moullin}).

\subsection{Higher rank graph $C^*$-algebras and canonical shift maps}

Whenever $\algen$ is a $C^*$-algebra, $\algen_h$ and  $\algen_+$ will denote  its hermitian and positive part respectively. For maps
acting between unital $C^*$-algebras the abbreviation ucp will mean unital and completely positive.

Given a higher rank graph $\Lambda$ the higher rank graph $C^*$-algebra $O_{\Lambda}$ is the universal  $C^*$-algebra generated by a family of
operators $s_u$, indexed by morphisms $u \in \Lambda$
satisfying
\begin{enumerate}
\item $\{ s_{a}=p_a \mid a \in B \}$ are pairwise orthogonal projections;
\item $s_{uv}= s_u s_v$ for all $u,v \in \Lambda$ such that $t(u)=o(v)$;
\item $s_u^* s_u = s_{t(u)}$ for all $u \in \Lambda$;
\item $s_a = \sum_{\sigma (u)=n, o(u)=a} s_u s_u^* $ for all $n \in \bz^r_+$.
\end{enumerate}
We will assume from now on that $\Lambda$ satisfies Kumjian and Pask's aperiodicity condition since in this case the algebra $O_{\Lambda}$ is nuclear and does not depend on the actual choice of generators (cf.\:\cite{kupa}). Under a mild additional assumption $O_{\Lambda}$ is also simple and purely infinite, as is stated in the same reference (the full proof of the fact that $O_{\Lambda}$ is purely infinite is given in \cite{Sims}).

There are two important subalgebras of $O_{\Lambda}$, the commutative $C^*$-algebra
$\com$ generated by the set $\{s_w s_w^*: w \in \Lambda \}$ and the AF $C^*$-algebra
$\hyp$ generated by the set $\{s_u s_w^*: u,w \in \Lambda , \sigma(u) = \sigma(w)\}$. The following fact is easy to check.

\begin{prop}
The algebra $\com$ is $*-$isomorphic to the algebra $C(X)=C(\Lambda^{\infty})$ of continuous functions on the
subshift $X$. The standard isomorphism is given by the linear extension of the map $s_u s_u^* \mapsto \chi_{Z_u}$.
\end{prop}

Whenever $f\in C(X)$, the corresponding element of $\com$ will be denoted by $\wt{f}$.
For further reference define (for $l \in \Zr$)
\[
\omega_l = \{s_u s_v^*: \sigma(u)\leq l \text{ and } \sigma(v) \leq l\}.
\]
For each $p \in \Zr$ define a noncommutative shift (canonical shift map) $\Phi^p:O_{\Lambda} \to O_{\Lambda}$
by
\[
\Phi^p (X) = \sum_{w\in \Lambda_p} s_w X s_w^*, \;\; X\in O_{\Lambda}.
\]
As the conditions above imply that $\sum_{w \in \Lambda_p} s_w s_w^* = 1_{O_{\Lambda}}$,
each $\Phi^p$ is a ucp map.
$\Phi^p$ defines an action of $\Zr$ on $O_{\Lambda}$ given by canonical shift maps.
Observe that the terminology is coherent with the shifts introduced in the previous section: for all $f\in C(X)$, $p \in \Zr$
\[
\Phi^p (\wt{f}) = \wt{f\circ T^p}.
\]

\section{Topological pressure and entropy for shifts on $O_{\Lambda}$}


Let $\algen$ be a  unital $C^*$-algebra.
We say that $(\phi,\psi,C)$ is an approximating triple for $O_{\Lambda}$ if $C$ is a
finite-dimensional $C^*$-algebra and both $\phi:C\to \algen$, $\psi:\algen \to C$ are
unital and completely positive. This will be indicated by writing $(\phi,\psi,C)\in CPA(\algen)$.
Whenever $\Omega$ is a finite subset of $\algen$ ($\Omega \in FS(\algen)$) and $\varepsilon >0$ the statement
 $(\phi,\psi,C)\in CPA(\algen, \Omega, \varepsilon)$ means that $(\phi,\psi,C)\in CPA(\algen)$ and
 for all $a \in \Omega$
 \[ \|\phi \circ \psi(a) - a \| < \varepsilon.\]
Nuclearity of $\algen$ is equivalent to the fact that for each $\Omega \in FS(\algen)$ and $\varepsilon >0$
there exists a triple  $(\phi,\psi,C)\in CPA(\algen, \Omega, \varepsilon)$. For such algebras  one can define
\[
\textrm{rcp}(\Omega, \varepsilon)= \min\{\textrm{rank}\,C:\; (\phi,\psi,C)\in CPA(\algen, \Omega, \varepsilon)\},
\]
where $\textrm{rank}\,C$ denotes the dimension of a maximal abelian subalgebra of $C$.
Let us recall the definition of topological pressure
in nuclear unital $C^*$-algebras, due to S.\,Neshveyev and E.\,St\o rmer (\cite{asympt}).
Assume that $\algen$ is nuclear, $a \in \algen_h$ and $\theta:\algen\to \algen$ is a ucp map.
For any $\Omega \in FS(\algen)$ and $n\in \bn$ let
\begin{equation}
\label{enning}
\text{orb}^n(\Omega)=\Omega^{(n)} = \bigcup_{j=0}^n {\theta^j(\Omega)},\;\;\; a^{(n)} = \sum_{j=0}^n {\theta^j(a)}.
\end{equation}
Define the \emph{noncommutative partition function}  ($\varepsilon>0, \, n\in\bn$)
\[
Z_{\theta,n} (a,\Omega,\varepsilon) = \inf \{\Tr \;\textrm{e}^{\psi(a^{(n)})}: \,
  (\phi,\psi,C)\in CPA(\algen, \Omega^{(n)}, \varepsilon)\},
\]
where $\Tr$ denotes a canonical trace on $C$ ($\Tr(q)=1$ for any minimal projection
$q\in C$).
Define
\[
P_{\theta} (a,\Omega,\varepsilon) =\limsup_{n\to \infty} \frac{1}{n}
\log \left( Z_{\theta,n} (a,\Omega,\varepsilon) \right),
\]
and the \emph{noncommutative pressure}
\[
P_{\theta} (a) = \sup_{\varepsilon>0, \, \Omega \in FS(\algen)}  P_{\theta} (a,\Omega,\varepsilon).
\]
If $a=0$ then $\Tr \;\textrm{e}^{\psi(a)}= \textrm{rank} \; C$ and  we recover Voiculescu's definition of
topological (approximation) entropy:
\[
\hte(\theta) = P_{\theta} (0) = \sup_{\varepsilon>0, \, \Omega \in FS(\algen)} \left(\limsup_{n\to\infty}
\left(\frac{1}{n} \log \textrm{rcp}(\Omega^{(n)}, \varepsilon)\right) \right).
\]
As shown in \cite{Voic} Proposition 4.8 the approximation entropy coincides with
classical topological entropy in the commutative case. The same holds for the pressure.
For this and extensions to the case of exact $C^*$-algebras and various related
topics we refer to \cite{Kerr}.

The computation of topological entropy and pressure for $\Phi^p$ hinges on the arguments
of \cite{Boca}, extended later in \cite{Kerr}. The essential fact is a possibility of embedding
$O_{\Lambda}$ in suitable matrix algebras of $O_{\Lambda}$
in such a way that the action of canonical shift maps takes a relatively simple form with respect
to these embeddings. (For Cuntz algebras this idea goes back to M.D.\,Choi).

Let $m\in \Zr$ and define $\rho_m:O_{\Lambda} \to M_{w_m} \ot O_{\Lambda}$ by
\[
\rho_m (X) = \sum_{u,w \in \Lambda_m} e_{u,w} \ot s_{u}^* X s_{w}, \;\;\; X \in O_{\Lambda},
\]
where $e_{u,w}$ denote standard matrix units in $M_{w_m}$. Elementary computations and uniqueness of
$O_{\Lambda}$ yield the following.
\begin{prop}
Each map $\rho_m$ is an injective $*$-homomorphism.
\end{prop}
The crucial fact is
\begin{lem} \label{main}
Let $u,w \in \Lambda$, $p \in \Zr$. Put $\maxshape = \sup \{\sigma(u),\sigma(w)\}$ (coordinatewise maximum). For each $m\in \Zr$ such that
$m\geq p+\maxshape$
\[
\rho_m (\Phi^p (s_u s_w^*)) = \sum_{\kappa \in \Lambda_{\maxshape - \sigma(w)}, \atop{ \lambda \in \Lambda_{\maxshape - \sigma(u)}}}
T_{\kappa,\lambda} \ot s_{\kappa} s_{ \lambda}^*,
\]
where $T_{\kappa,\lambda} \in M_{w_m}$ are partial isometries.
\end{lem}

\begin{proof}
In the following we use the convention $s_{\mu \nu}:=0$ if $\mu, \nu \in \Lambda$ and $t(\mu) \neq o(\nu)$.
Note that
\[
\Phi^l (s_u s_w^*) = \sum_{\nu \in \Lambda_l} s_{\nu} s_{u}s_{w}^* s_{\nu}^* =
\sum_{\nu \in \Lambda_l} s_{\nu u} s_{\nu w}^*,
\]
and further, using the defining relations of $O_{\Lambda}$
\begin{eqnarray*}
\rho_m (\Phi^l (s_u s_w^*))
&=&
\sum_{a,b \in \Lambda_m} e_{a,b} \ot \sum_{\nu \in \Lambda_l} s_a^* s_{\nu u} s_{\nu w}^* s_{b} \\
&=&
\sum_{a,b \in \Lambda_m,\, \nu \in \Lambda_l,\atop{\gamma\in \Lambda_{m-l-\sigma(u)}}}
{e_{a,b} \ot s_{a}^* s_{\nu u \gamma} s_{\nu w \gamma}^* s_{b}} \\
&=&
\sum_{b \in \Lambda_m,\, \nu \in \Lambda_l,\atop{\gamma\in \Lambda_{m-l-\sigma(u)}}}
 {e_{\nu u \gamma,b} \ot s_{\nu w \gamma}^* s_{b}} \\
&=&
\sum_{b \in \Lambda_m,\, \nu \in \Lambda_l,\, \gamma \in \Lambda_{m-l-\sigma(u)},\atop{\kappa \in \Lambda_{\maxshape-\sigma(w)}, \, \delta \in \Lambda_{\maxshape-\sigma(u)}}}
{e_{\nu u \gamma,b} \ot s_{\kappa} s_{\nu w \gamma\kappa}^* s_{b\delta} s_{\delta}^*} \\
&=&
\sum_{\nu \in \Lambda_l,\, \gamma \in \Lambda_{m-l-\sigma(u)}, \atop{\kappa \in \Lambda_{\maxshape-\sigma(w)}}}
{e_{\nu u \gamma,\nu w \gamma \kappa|_{m]}} \ot s_{ \kappa} s^*_{\nu w \gamma\kappa|_{[m}}} \\
&=&
\sum_{\kappa \in \Lambda_{\maxshape- \sigma(w)}, \atop{ \lambda \in \Lambda_{\maxshape -\sigma(u)} }}
\left(\sum_{\nu \in \Lambda_l,\, \gamma\in \Lambda_{m-l-\sigma(u)},
\atop{ \nu w \gamma\kappa|_{[m}= \lambda }}
e_{\nu u \gamma,\nu w \gamma\kappa|_{m]}} \right) \ot s_{\kappa} s^*_{\lambda}.
\end{eqnarray*}
Define now for $\kappa \in \Lambda_{\maxshape - \sigma(w)}, \,
\lambda \in \Lambda_{\maxshape - \sigma(u)}$
\[
T_{\kappa,\lambda} =
\sum_{\nu \in \Lambda_l,\, \gamma\in \Lambda_{m-l-\sigma(u)},
\atop{ \nu w \gamma\kappa|_{[m} = \lambda }}
e_{\nu u \gamma,\nu w \gamma\kappa|_{m]}}.
\]
It remains to check that the matrices thus defined are indeed partial isometries. This, however, may be seen by exploiting the fact
they are defined in terms of matrix units (when represented in a basis indexed by words in $\Lambda_m$ the matrices
$T_{\kappa,\lambda}$ have at most one nonzero element, equal to $1$, in each row and column).
\end{proof}

Before we formulate the main theorem of this section we need a multidimensional version of Proposition 5.1 in \cite{Kerr}:
\begin{prop} \label{ineq}
Let $\wt{f} \in \com_+$, $u,w \in W$ with $\sigma(u)=\sigma(w)$. If $u \neq w$ then $s_{u}^* \wt{f} s_w =0.$ Moreover
\[ s_{u}^* \wt{f} s_{u} \leq \max \{f(x):x \in Z_u\} 1_{O_{\Lambda}}.\]
\end{prop}
\begin{proof}
This follows since $s_us_u^* \wt{f} s_w s_w^*=  \wt{f} s_us_u^* s_w s_w^*$ and thus $s_us_u^* \wt{f} s_u s_u^*   \leq \max \{f(x):x \in Z_u\} s_u s_u^* \leq \max \{f(x):x \in Z_u\} 1_{O_{\Lambda}}$.
\end{proof}
\begin{tw}
Let $\wt{f} \in \com_h$. Then
\[
P_{\Phi^p} (\wt{f}) = P_{\textrm{top}} (f,T^p).
\]
for every $p \in \Zr$.
\end{tw}

\begin{proof} The proof is an adaptation of arguments of \cite{Boca} and \cite{Kerr}. We may and do assume that $f\geq 0$.
The monotonicity of pressure proved in \cite{Kerr} yields
\[
P_{\textrm{top}} (f,T^p) \leq P_{\Phi^p} (\wt{f}).
\]
To obtain the reverse inequality fix $l\in \Zr$, $\delta >0$. As $O_{\Lambda}$  is nuclear, there exists a triple
$(\phi_0,\psi_0, M_{C_l}) \in CPA (O_{\Lambda},\omega_l, \frac{1}{4w_{l}}\delta)$.
Fix $n\in \Zr$, put $m=n+l$ and consider the following diagram (compare with \cite{Boca})

\begin{picture}(350,100)
\put(0,80){$O_{\Lambda}$}
\put(80,80){$\rho_m(O_{\Lambda})$}
\put(92,50){$M_{w_m} \otimes O_{\Lambda}$}
\put(325,80){$O_{\Lambda}$}
\put(205,80){$\rho_m(O_{\Lambda})$}

\put(193,50){$M_{w_m} \otimes O_{\Lambda}$}

\put(262,64){$M_d$}

\put(148,0){$M_{w_m} \otimes M_{C_l}$}

\put(39,88){$\rho_m$}
\put(20,83){\vector(1,0){50}}

\put(96,73){\vector(1,-1){14}}

\put(129,32){$\textrm{id} \ot \psi_0$}
\put(116,42){\vector(1,-1){27}}

\put(176,32){$\textrm{id} \ot \phi_0$}
\put(192,15){\vector(1,1){27}}

\put(224,65){\line(0,1){5}}
\put(216,65){$\cap$}

\put(245,78){\vector(2,-1){14}}
\put(245,68){$\gamma$}

\put(280,69){\vector(3,1){30}}
\put(292,66){$\eta$}

\put(271,89){$\rho_m^{-1}$}
\put(247,84){\vector(1,0){62}}

\put(243,55){\vector(2,1){16}}
\put(249,51){$\mu$}

\put(17,71){\vector(2,-1){123}}
\put(65,33){$\psi$}

\put(205,11){\vector(2,1){113}}
\put(265,31){$\phi$}

\end{picture}

\vspace{0.4cm}

The embedding $\rho_m(O_{\Lambda}) \hookrightarrow M_{w_m} \otimes O_{\Lambda}$ does not allow a ucp left-inverse but approximate
inverses (using Arveson's extension Theorem) as indicated in the diagram.
More precisely, nuclearity of $O_{\Lambda}$ implies that there exists $d \in \bn$ and ucp maps
$\gamma:\rho_m(O_{\Lambda})\to M_d$ and $\eta:M_d\to O_{\Lambda}$ such that for all $a\in \rho_m(\Phi^k(\omega_l))$,
$k\leq n$,
\[
\| \eta \circ \gamma (a) - \rho_m^{-1} (a) \| < \frac{\delta}{2}.
\]
Let $\mu:M_{w_m} \ot O_{\Lambda} \to M_d$ be a ucp extension of $\gamma$.
Consider again  any $X \in \omega_l$ and let $k\in \Zr$, $k\leq n$. Then
\begin{eqnarray*}
\lefteqn{\|\phi \circ \psi  (\Phi^k(X)) - \Phi^k(X) \|}  \hspace*{1.5cm}
 \\
&=&\|\eta \circ \mu \circ (\textrm{id} \ot \phi_0 \circ \psi_0)\circ \rho_m  (\Phi^k(X))
- (\rho_m^{-1} \circ \rho_m)(\Phi^k(X)) \| \\
& \leq & \|\eta \circ \mu \circ (\textrm{id} \ot \phi_0 \circ \psi_0)\circ \rho_m  (\Phi^k(X)) - \eta \circ \mu \circ \rho_m (\Phi^k(X))\| \\
& & \;\;\;\;\;\;\;\;\;\;\;\;\;\;\;\;\;\;\;\;\;\;
+ \;\;\|\eta \circ \mu \circ \rho_m (\Phi^k(X)) - (\rho_m^{-1} \circ \rho_m)(\Phi^k(X))\| \\
&< &\|  (\textrm{id} \ot \phi_0 \circ \psi_0) \circ \rho_m  (\Phi^k(X))
- \rho_m(\Phi^k(X)) \| +\frac{\delta}{2}.
\end{eqnarray*}
As $X=s_u s_w^*$ for some $u,w \in \Lambda, \sigma:=\max\{\sigma(u), \sigma(w)\} \leq l$, we are in a position to apply Lemma \ref{main} to obtain
\begin{eqnarray*}
\lefteqn{\|\phi \circ \psi  (\Phi^k(X)) - \Phi^k(X) \| } \hspace*{1.5cm}\\
&< & \Big\| \sum_{\kappa \in \Lambda_{\sigma - \sigma(w)}, \atop{ \lambda \in \Lambda_{\sigma - \sigma(u)}}}
T_{\kappa,\lambda} \ot \big( (\phi_0 \circ \psi_0) ( s_{\kappa}s_{\lambda}^*) - s_{\kappa} s_{\lambda}^* \big)  \Big\|+
 \frac{\delta}{2} \\
& < &
2 w_l \frac{\delta}{4w_l}+ \frac{\delta}{2} = \delta,
\end{eqnarray*}
and we proved that
\begin{equation}
\label{approx} (\phi,\psi,M_{w_m} \ot M_{C_l}) \in CPA (O_{\Lambda}, \omega_l^{(n)}, \delta),
\end{equation}
where $\omega_l^{(n)}=\bigcup_{k\in\Zr:k\leq n} \Phi^k(\omega_l)$.

Let now $s \in \bn$ and, following the convention in (\ref{enning}) for $\theta=\Phi^p$, put $\wt{f}^{(s)} = \sum_{t=0}^s \Phi^{tp} (\wt{f})$.
Let $\text{orb}^s (\omega_l)=\bigcup_{t=0}^s \Phi^{tp}(\omega_l)$. As $\text{orb}^s (\omega_l) \subset \omega_l^{(sp)}$, the statement (\ref{approx}) for $n=sp$
implies  $(\phi_m,\psi_m, M_{w_m} \ot M_{C_l}) \in CPA (O_{\Lambda}, \text{orb}^s (\omega_l) , \delta)$. Moreover
\begin{eqnarray*}
\textrm{Tr} \exp \left(\psi (\wt{f}^{(s)})\right)
&=&
\textrm{Tr} \exp \left((\textrm{id} \ot \psi_0) (\rho_m (\wt{f}^{(s)}))\right)  \\
&=&
\sum_{q=0}^{\infty} \frac{1}{q!} \textrm{Tr} \left( \left(\sum_{u,w \in \Lambda_m} e_{u,w} \ot
\psi_0 (s_{u}^* \wt{f}^{(s)} s_w)\right)^q \right)
\end{eqnarray*}
and by Proposition \ref{ineq}
\begin{eqnarray*}
\lefteqn{ \textrm{Tr} \exp \left(\psi (\wt{f}^{(s)})\right)   } \\
&=&
\sum_{u \in \Lambda_m} \textrm{Tr} \exp \left(e_{u,u} \ot \psi_0 (s_{u}^* \wt{f}^{(s)} s_u )\right) \\
&\leq&
\sum_{u \in \Lambda_{sp+l}} \textrm{Tr} \exp \left( \big(e_{u,u} \ot 1_{M_{C_l}} \big) \max\{ f^{(s)}(x) :   x \in Z_u \} \right) \\
&=&
C_l
\sum_{u \in \Lambda_{sp+l}} \exp\left(\max\{ f(x) +f(T^p(x)) + \ld +f(T^{sp}(x)): x \in Z_u \}\right).
\end{eqnarray*}
The fact that $\bigcup_{l\in \Zr} \bigcup_{q=0}^{\infty} {\Phi^{qp}(\omega_l)}$ is total in $O_{\Lambda}$,
the Kolmogorov-Sinai property for noncommutative pressure (Proposition 3.4 of \cite{Kerr}) and Remark \ref{compre} end the proof.
\end{proof}

\noindent
The above theorem (put $f=0$) and Proposition \ref{coment} yield immediately:
\begin{cor}
For a higher rank Cuntz-Krieger algebra the topological entropy of $\Phi^p$ is equal to $\log \textrm{r}(M^p)$. Moreover,
the topological entropy of $\Phi_j:=\Phi^{e_j}$ is equal to $\log \textrm{r}(M_j)$ for every $j \in \{1,\ldots ,r \}$.
\end{cor}
\begin{rem}
When $r=2$ and  $M_1,M_2$ are of the form $M_1 = A_1 \ot I$, $M_2=I \ot A_2$ for some matrices $A_1\in M_{d_1}$, $A_2\in M_{d_2}$,
it is well known that the rank 2 Cuntz-Krieger algebra is given by $O_M=O_{A_1}  \ot O_{A_2}$.
Then $M_1 M_2=A_1 \ot A_2$, and $\log \textrm{r}(M_1 M_2)= \log \textrm{r}(M_1) + \log \textrm{r}(M_2)$.
It would be interesting to find an example of matrices $M_1,M_2$ satisfying (H0)-(H3) such that
$\log \textrm{r}(M_1M_2)\neq \log \textrm{r}(M_1) + \log \textrm{r}(M_2)$; this being equivalent to
$\hte (\Phi_{\textup{diag}}) \neq \hte(\Phi_1) + \hte(\Phi_2)$, with $\Phi_{\textup{diag}}$ denoting the diagonal
shift. Note that we always have
$\log \textrm{r}(M_1M_2) \leq \log \textrm{r}(M_1)+ \log \textrm{r}(M_2)$ as follows from elementary
properties of the spectral radius.
\end{rem}

\section{Entropy of the $\Zr$-action on $O_{\Lambda}$}
Extensions of a classical concept of topological entropy to actions of amenable groups date back to 1970s.
Here we are mainly interested in the special case of $\Zr$ actions.
Let $\algen$ be a $C^*$-algebra and assume that $\alpha$ is an action of $\Zr$ on $\algen$ by unital completely positive maps.
For each $\Omega \in FS(\algen)$ and $n\in \bn$ put $\Omega^{(n)} = \bigcup_{l\in \Zr, l\leq \bar{n}}
\alpha^l(\Omega)$ and
\[
\hte (\alpha,\Omega, \varepsilon) = \limsup_{n\to \infty} \left(\frac{1}{n^r} \log \textrm{rcp}(\Omega^{(n)}, \varepsilon)\right).
\]
We define the \emph{topological entropy of the action $\alpha$} by
\[
\hte (\alpha) =\sup_{\Omega\in FS(\algen), \varepsilon>0} \hte (\alpha,\Omega, \varepsilon).
\]
This definition coincides with the standard one in the commutative case:
\begin{prop}
Let $\alpha$ be an action of $\Zr$ on a compact metric space $X$  by continuous transformations.
Let $\wt{\alpha}$ be the corresponding action of
$\Zr$ on $C(X)$:
\[ \wt{\alpha^l} (f) = f\circ \alpha^l ,\;\;\; f \in C(X),\, l\in \Zr.\]
Then
\[ h_{\textrm{top}} (\alpha) = \hte(\wt{\alpha}).\]
\end{prop}

The above proposition may be proved along the same lines as Proposition 4.8 in \cite{Voic} - the arguments there apply also
for continuous transformations which are not neccessarily homeomorphic, and the classical variational principle holds for
$\Zr$-actions (see e.g.\
\cite{Misi}).
The other important fact concerning the entropy defined as above is the Kolmogorov-Sinai property, which again can be proved by a
straightforward modifications of arguments in \cite{Voic}:

\begin{prop}\label{Kolsin}
Let $\{\Omega_k: k \in \bn\}$ be an increasing sequence of finite subsets of $\algen$ such that $\bigcup_{k\in \bn} \bigcup_{l\in \Zr}
{\alpha^l (\Omega_k)}$ is total in $\algen$. Then
\[
\hte (\alpha) =\sup_{k \in \bn, \varepsilon>0} \hte (\alpha,\Omega_k , \varepsilon).
\]
\end{prop}
With that at hand, we can show the following:
\begin{tw} \label{someth}
Let $r\geq 2$ and $O_{\Lambda}$ a rank-r graph $C^*$-algebra. The topological entropy of the action $\alpha$ on $O_{\Lambda}$ by shifts is $0$.
\end{tw}

\begin{proof}
Let $k\in \bn$ and consider the set $\omega_{\bar{k}}$.  Fix $\delta>0$. By Proposition \ref{Kolsin} it is enough to prove that
\begin{equation}\hte (\alpha,\omega_{\bar{k}}, \delta) = 0. \label{zero}\end{equation}
However, (\ref{approx}) yields (for all $n \in \bn$)
\[ \textrm{rcp}(\omega_{\bar{k}}^{(n)}, \delta) \leq C_{\bar{k}} w_{\overline{n+k}},\]
the factorisation property implies that
\[ w_{{n+k}} \leq w_{(n+k)e_1} \cdots w_{(n+k)e_r} \leq (\max_{i=1,\ldots,r} {|\Lambda_{e_i}|} ) ^{(n+k)r}\]
and (\ref{zero}) follows.
\end{proof}


\end{document}